\documentclass{article}
\usepackage[utf8]{inputenc}
\usepackage{amsmath,amssymb,amsthm}
\usepackage{enumerate}
\usepackage{url}
\newtheorem{thm}{Theorem}[section]
\newtheorem{dfn}[thm]{Definition}
\newtheorem{lem}[thm]{Lemma}
\newtheorem{prop}[thm]{Proposition}

\theoremstyle{remark}
\newtheorem{rem}[thm]{Remark}
\newtheorem{ex}[thm]{Example}

\DeclareMathOperator{\scal}{Scal}
\DeclareMathOperator{\ric}{Ric}
\DeclareMathOperator{\R}{R}
\DeclareMathOperator{\Lop}{L}
\DeclareMathOperator{\diam}{diam}
\DeclareMathOperator{\vol}{vol}
\DeclareMathOperator{\inj}{inj}

\DeclareMathOperator{\tr}{trace}
\DeclareMathOperator{\Id}{I}
\DeclareMathOperator{\id}{id}

\newcommand{\ddt}[1]{\frac{\partial #1 }{\partial t}}
\newcommand{\ps}[2]{\left\langle #1,#2 \right\rangle}
\newcommand{\eps}{\varepsilon}
\newcommand{\ph}{\varphi}

\author{Thomas Richard}
\title{Lower bounds on Ricci flow invariant curvatures and geometric applications.}
\begin{document}
\maketitle
\begin{abstract}
We consider Ricci flow invariant cones $\mathcal{C}$ in the space of curvature
operators lying between nonnegative Ricci curvature and nonnegative
curvature operator. Assuming some mild control on the scalar curvature
of the Ricci flow, we show that if a solution to Ricci flow has its
curvature operator which satsisfies $\R+\eps \Id\in\mathcal{C}$ at the
initial time, 
then it satisfies $\R+K\eps\Id\in\mathcal{C}$ on some time interval
depending only on 
the scalar curvature control. 

This allows us to link Gromov-Hausdorff
convergence and Ricci flow convergence when the limit is smooth and
$\R+\Id\in\mathcal{C}$ along the sequence of initial
conditions. Another application is a stability result for 
manifolds whose curvature operator is almost in $\mathcal{C}$.

Finally, we study the case where $\mathcal{C}$ is contained in the
cone of operators whose sectional curvature is nonnegative. This allow
us to weaken the assumptions of the previously mentioned
applications. In particular, we construct a Ricci flow for a class of
(not too) singular Alexandrov spaces.
\end{abstract}
% MSC2010 Classification : 53C44, 53C20.
\section{Introduction and statement of the results}

In the study of the Ricci flow, various nonnegative curvature
conditions have been shown to be preserved, and the discovery of new
invariant conditions has often given rise to new geometric
applications. One of the most famous occurence of this fact is the
discovery by Brendle and Schoen in \cite{MR2449060} and independently by
Nguyen in \cite{MR2587576} of the preservation of
nonnegative isotropic curvature, which plays a crucial role in the
proof by Brendle and Schoen of the differentiable sphere theorem in \cite{MR2449060}.

Once one has understood the behaviour of the Ricci flow assuming the
nonnegativity of a certain curvature, it is natural to ask if
something can be done under arbitrary lower bounds on this given
curvature. Such a work has been done for Ricci curvature in 
dimension 3 by Simon in \cite{MR2526789} and \cite{MS2009}. An
important feature of this work
is that, in order to control lower bounds on the Ricci curvature along
the flow, one has to impose further geometric conditions on the
initial manifold. In Simon's work, a non-collapsing assumption is
required. Our estimate will rely on an a priori bound on the scalar
curvature. 

In order to state the results of this paper, we need some
terminology.
\begin{dfn}
  A nonnegativity condition on the curvature is given by a closed
  convex cone $\mathcal{C}$ in the space of algebraic curvature
  operators $S^2_B\Lambda^2\mathbb{R}^n$ such that :
  \begin{itemize}
  \item The identity operator
    $\Id:\Lambda^2\mathbb{R}^n\to\Lambda^2\mathbb{R}^n$ lies in the
    interior of $\mathcal{C}$.
  \item $\mathcal{C}$ is invariant under the action of
    $O(n,\mathbb{R})$ on $S^2_B\Lambda^2\mathbb{R}^n$ given by :
    \[ \ps{g.\R(x\wedge y)}{z\wedge t}=\ps{\R(gx\wedge gy)}{gz\wedge
      gt}.\]
  \end{itemize}
\end{dfn}
Recalls and references on algebraic curvtaure operators are included
in Section \ref{sec:prel-about-algebr}.

Given a nonnegativity condition $\mathcal{C}$ and a Riemannian
manifold $(M,g)$, we can canonically embed $\mathcal{C}$ in
$S^2_B\Lambda^2T_mM$ for each $m\in M$, thanks to the
$O(n,\mathbb{R})$ invariance of $\mathcal{C}$. We say that $(M,g)$ has
$\mathcal{C}$-nonnegative curvature (or $\R\geq_\mathcal{C}0$) if, for
each $x\in M$ the curvature operator of $(M,g)$ at $x$ belongs to
$\mathcal{C}$. Classical condtions of nonnegative curvature operator,
nonnegative sectional curvature, nonnegative Ricci curvature or
nonnegative scalar curvature fit in this framework.

Similarly we say that $(M,g)$ has $\mathcal{C}$-curvature bounded from
below by $-k\Id$ (or $\R\geq_\mathcal{C} -k\Id$) for some
$k\in\mathbb{R}$ if for each $x\in M$ the
curvature operator $\R$ at $x$ is such that $\R+k\Id\in\mathcal{C}$.

We now define a class of nonnegativity condition which behaves well
with Ricci flow.
\begin{dfn}
  A nonnegativity condition $\mathcal{C}$ is said to be (Ricci Flow)
  invariant if $\mathcal{C}$ is preserved by Hamilton's ODE
  $\dot{\R}=2Q(\R)$. Namely, if $\R(t)$ is a solution to Hamilton's ODE
  on some time interval $[0,T)$ such that $\R(0)\in\mathcal{C}$, then
  $\R(t)\in\mathcal{C}$ for all $t\in[0,T)$.
\end{dfn}
Details and references about Hamilton's ODE are given in Section
\ref{sec:prel-about-algebr}.

Hamilton's maximum principle for tensors (\cite[Theorem
4.3]{MR862046}) implies that such a cone is 
preserved by Ricci flow in the sense that, if $(M,g_0)$ is a compact
Riemannian manifold such that $\R\geq_\mathcal{C} 0$, then the Ricci
flow $(M,g(t))$ such that $g(0)=g_0$ satisfies
$\R(g(t))\geq_\mathcal{C} 0$ as long as it exists.

We are now ready to state our result. It roughly says the
following. We consider a manifold whose
$\mathcal{C}$-curvature is bounded from below, where $\mathcal{C}$ is an
invariant condition between nonnegative Ricci curvature and
nonnegative curvature operator. We furthermore assume that an a priori estimate
on the blow up rate of the scalar curvature of the Ricci flow as $t$
goes to zero is true. Then the $\mathcal{C}$-curvature can be bounded from
below on a small time interval.
\begin{thm}\label{lboundthm}
  For any dimension $n\in\mathbb{N}$, any $A\in(0,\frac{1}{4})$
  and any $B>0$, one can find $T=T(n,A,B)$ and $K=K(n,A,B)$
  such that if 
  $\mathcal{C}\subset S^2_B\Lambda^2\mathbb{R}^n$ is a closed convex
  cone which satisfies : 
    \begin{enumerate}
    \item $\mathcal{C}$ is an invariant nonnegativity condition,
    \item $\mathcal{C}$ contains the cone of nonnegative curvature operators,
    \item $\mathcal{C}$ is contained in the cone of curvature operators 
      whose Ricci curvature is nonnegative,
    \end{enumerate}
  and $(M^n,g(t))_{t\in [0,T')}$ is a Ricci flow on a smooth compact
  manifold satisfying :
    \begin{enumerate}
    \item $\R(g(0))\geq_{\mathcal{C}}-\eps \Id$ at $t=0$ 
      for some $\eps\in [0,1]$,
    \item $|\scal(g(t))|\leq A/t+B$ for $t$ in $(0,T')$,
    \end{enumerate}
  we have :
\[\R(g(t))\geq_{\mathcal{C}}-K\eps \Id\]
  for all $t$ in $[0,T')\cap [0,T)$.
\end{thm}
\begin{rem}
  During the redaction of this article, the author has been informed
  that a similar estimate was also known by Miles Simon and Arthur Schlichting.
\end{rem}
\begin{ex}\label{excones}
  Known examples of cones which satisfy the assumptions of the theorem
  include :
  \begin{itemize}
  \item the cone $\mathcal{C}_{CO}$ of nonnegative curvature operators,
  \item the cone $\mathcal{C}_{2CO}$ of 2-nonnegative curvature operators,
  \item the cone $\mathcal{C}_{IC1}$ of curvature operators which
    have nonnegative isotropic when extended by $0$ to $\Lambda^2\mathbb{R}^{n+1}$,
  \item the cone $\mathcal{C}_{IC2}$ of curvature operators which
    have nonnegative isotropic when extended by $0$ to $\Lambda^2\mathbb{R}^{n+2}$.
  \end{itemize}
  All these conditions have been extensively studied 
  (\cite{MR2415394},\cite{MR2449060},\cite{MR2386107},\cite{MR2462114}) 
  and compact manifolds with $\mathcal{C}$-nonnegative
  curvature have been classified when $\mathcal{C}$ is one of these
  four cones. An exposition of the relations
  between these 
  condtions and how nonnegativity of these curvatures affect the
  topology of the underlying manifold can be found in the Brendle's
  book \cite{MR2583938} together with precise definitions and
  additional references. It should also be noted that Wilking has
  given a unified proof of the preservation of these conditions (along
  with others) in \cite{2010arXiv1011.3561W}.

  Some continuous families of such cones have also been constructed in
  \cite{MR2415394} and \cite{2011arXiv1105.5311G}.

  It should be noted that in dimension greater or equal to $4$, nonnegative Ricci
  curvature is not preserved, see \cite{MR2736347}.
\end{ex}
\begin{rem}\label{remwilk}
  If $\mathcal{C}$ satisfies the assumptions of the theorem and
  moreover is a Wilking cone (see \cite{2010arXiv1011.3561W}), it follows
  from the work of Gururaja, Maity and Seshadri in
  \cite{2011arXiv1101.5884G} that 
  $\mathcal{C}$ is included in $\mathcal{C}_{IC1}$.
\end{rem}

{\bf From now on any curvature condition $\mathcal{C}$ is supposed to satisfy
the assumptions of Theorem \ref{lboundthm}.}

The estimate of Theorem \ref{lboundthm} allows us to adapt the methods
of \cite{MR2526789} and \cite{MS2009} to some
higher dimensional situations. 

In the first two applications, the
estimate on the scalar curvature which is required to apply Theorem
\ref{lboundthm} will be obtained by Perelman's pseudolocality theorem,
first stated in \cite{PerelEnt} but omitting the crucial assumption
of completeness as pointed by Topping, see \cite[Theorem
A.3]{2011arXiv1106.2493G}, 
complete statement and proofs can be found in
\cite{KleiLott,MR2779131}.

Our first application is to show that,
if the $\mathcal{C}$-curvature is bounded from below along a
sequence of compact $n$-dimensional smooth manifolds which Gromov-Hausdorff
converges (we will write GH-converges in the sequel) to a compact
$n$-dimensional smooth manifold, then, up to a subsequence, the
associated Ricci flows converge to a Ricci flow of the limit manifold
(where the initial condition is to be understood in a weak
sense). Here convergence is smooth convergence of the Ricci flows up
to diffeomorphisms, as in \cite{HamComp}. More precisely, we prove
the following theorem, which is an 
higher dimensional analogue of \cite[Theorem 9.2]{MR2526789}, where
such a theorem has been proved in dimension 3 under lower bounds on
the Ricci curvature and without assuming smoothness and compactness of
the limit :
\begin{thm}\label{GHrfthm}
  Let $(M_k,g_k)$ be a sequence of compact $n$-manifolds which 
  satisfies $R\geq_\mathcal{C}-\Id$ and which GH-converges
  to a compact smooth $n$-manifold $(M,g)$. Let
  $(M_k,g_k(t))_{t\in[0,T_k)}$ be the maximal solution of the Ricci
  flow satisfying $g_k(0)=g_k$. Then :
  \begin{enumerate}
  \item there is a positive constant $T$ such that each Ricci flow
    $(M_k,g_k(t))$ is defined at least on $[0,T)$ and the sequence of
    Ricci flows $(M_k,g_k(t))_{t\in(0,T)}$ is precompact in the sense
    of Cheeger-Gromov-Hamilton.
  \item any Cheeger-Gromov-Hamilton limit
    $(\tilde{M},\tilde{g}(t))_{t\in(0,T)}$ of a convergent subsequence
    of $(M_k,g_k(t))_{t\in(0,T)}$ is such that $\tilde{M}$ is
    homeomorphic to $M$ and the distance functions $d_{\tilde{g}(t)}$
    uniformly converge as $t$ goes to $0$ to some distance $\tilde{d}$
    which is isometric to the distance $d_g$. In particular the
    $M_k$'s are homeomorphic to $M$ for $k$ large enough.
  \end{enumerate}
\end{thm}
\begin{rem}
  Along the proof ot Theorem \ref{GHrfthm}, we will see that the
  precompactness of the sequence of flows $(M_k,g_k(t))_{t\in(0,T)}$
  still holds when one replaces $\mathcal{C}$-curvature bounded from
  below by Ricci curvature bounded from below (see Lemma
  \ref{GHvsRF}). However, our method of proof requires the lower
  bound on the $\mathcal{C}$-curvature to control
  the initial condition of the limit flow.
\end{rem}
\begin{rem}
  In the conclusions of the theorem, the fact that the $M_k$'s are
  homeomorphic to $M$ for $k$ large enough can be seen using Cheeger
  and Colding's work on manifolds with Ricci curvature bounded from
  below (see \cite[Theorem A.1.12]{CheeCol}). Additionaly, Cheeger and
  Colding's result allow to strengthen the conclusion from
  homeomorphism to diffeomorphism. However, our proof is
  independent of this work.
\end{rem}

Another application is a result about manifolds whith almost
nonnegative $\mathcal{C}$-curvature, in the spirit of \cite[Theorem 1.7]{MR2526789} :

\begin{thm}\label{almostNNCC}
  For any $i>0$ and $D>0$, for any $n\in\mathbb{N}$, there is an $\varepsilon>0$ such that any manifold $(M,g)$ satisfying :
  \begin{enumerate}
  \item $\inj(g)\geq i$
  \item $\diam(M,g)\leq D$
  \item $\R\geq_{\mathcal{C}}-\varepsilon \Id$
  \end{enumerate}
admits a metric whose curvature is $\mathcal{C}$-nonnegative.
\end{thm}
Using the classifcation results of Brendle \cite[Theorem 9.33]{MR2583938}, Micallef
and Wang \cite[Theorem 3.1]{MR1253619} and remark \ref{remwilk}, if we moreover 
assume that $\mathcal{C}$ is a Wilking
cone, we have that the universal cover of $M$ is diffeomorphic to a
product $\mathbb{R}^k\times N_1\times\dots\times N_l$ where each $N_i$
is one of the following :
\begin{itemize}
\item a standard sphere $\mathbb{S}^n$ with $n\geq 2$,
\item a compact symmetric space.
\end{itemize}

We then impose stronger requirements on the cone $\mathcal{C}$. We
assume that $\mathcal{C}$ is included in the cone of curvature
operators whose sectional curvature is nonnegative. 
The cones which satisfy this assumptions in the list of Example
\ref{excones} are $\mathcal{C}_{CO}$ and $\mathcal{C}_{IC2}$.

This allows us to
weaken the hypothesis of our results. In this context, it turns out
that a convenient assumption that can be used to fulfill the
hypothesis of Theorem \ref{lboundthm} is that balls have almost
euclidean volume. This is proved in Lemma \ref{lemalmosteuc}, and was
inspired to the author by the recent work of Cabezas-Rivas and Wilking
\cite{2011arXiv1107.0606C}. For instance, Theorem \ref{GHrfthm} becomes :
\begin{thm}\label{GHRFnnsec}
  Let $\mathcal{C}$ be a cone satifying the hypothesis of Theorem
  \ref{lboundthm} and which is contained in the cone of curvature
  operator whose sectional curvature is nonnegative.
  
  For any $n\in\mathbb{N}$, there exist $\kappa>0$, $T>0$ and $\delta>0$
  such that if $(X,d)$ is a 
  metric space which is a Gromov-Hausdorff limit of a sequence of
  compact manifolds $(M_i^n,g_i)$ such that :
  \begin{itemize}
  \item $\R(g_i)\geq_{\mathcal{C}} -\kappa\Id$,
  \item for any $x\in M_i^n$, $\vol_{g_i}(B_{g_i}(x,1))\geq
    (1-\delta)\omega_n$, where $\omega_n$ is the volume of the unit
    ball in $\mathbb{R}^n$,
  \end{itemize}
  then one can find a Ricci flow $(M,g(t))$ defined on $(0,T)$ with
  bounded curvature on each time slice such
  that $M$ is homeomorphic to $X$ and the distance $d_g(t)$ converge
  uniformly on any compact of $M$ to a distance $\tilde{d}$ such that
  $(M,\tilde{d})$ is isometric to $(X,d)$. 
\end{thm}
\begin{rem}
  The fact that $X$ is a manifold is a direct consequence of
  Perelman's stability theorem (see \cite{MR2408265}), but we will not
  use this result in the proof. The metric
  space $(X,d)$ in our result is an Alexandrov space with curvature
  bounded form below and can have cone-like singularities, but the
  almost euclidean volume condition forbids
  too sharp cone angles.
\end{rem}

Similarly, we get a stronger analogue of Theorem \ref{almostNNCC} :
\begin{thm}\label{almostnnccnnsc}
  Let $\mathcal{C}$ be a cone satifying the hypothesis of Theorem
  \ref{lboundthm} and which is contained in the cone of curvature
  operator whose sectional curvature is nonnegative.

  For any $n\in\mathbb{N}$, there exists $\delta>0$ such that for any
  $D>0$, one can find $\eps>0$ such that if $(M^n,g)$ is a compact
  Riemannian manifold such that :
  \begin{itemize}
  \item $\R(g)\geq_\mathcal{C}-\eps\Id$,
  \item $\forall x\in M,\ \vol_g(B_g(x,1))\geq(1-\delta)\omega_n$,
  \item $\diam(M,g)\leq D$,
  \end{itemize}
  then $M$ admits a metric with $\mathcal{C}$-nonnegative curvature.
\end{thm}
The article is organised as follows : in Section
\ref{sec:prel-about-algebr} we set up the notations and give some
background about the evolution equation of the curvature operator
along the Ricci flow that will be used in the proof of Theorem
\ref{lboundthm}. In Section \ref{sec:proofthm}, we give the proof of
Theorem \ref{lboundthm}. The applications are discussed in Section
\ref{sec:applications}.  Section \ref{sec:manif-with-mathc} is devoted
to the applications in the case where $\mathcal{C}$-nonnegative
curvature implies nonneagtive sectional curvature.

\subsection*{Acknowledgements}

The author is grateful to Gilles Carron and Harish Seshadri for
helpful discussions during the elaboration of this paper. The author
also thanks his supervisor Gérard Besson for his interest and support.

\section{Preliminaries about algebraic curvature operators.}
\label{sec:prel-about-algebr}

In this section, we set up the notations that will be used in this
paper. Our conventions follow closely those of Böhm and Wilking
in \cite{MR2415394}. 

We will denote by $S^2\Lambda^2\mathbb{R}^n$ the vector space of
symmetric operators on $\Lambda^2\mathbb{R}^n$ equiped with the
standard inner product. $S^2_B\Lambda^2\mathbb{R}^n$ is the vector
space of operators in $S^2\Lambda^2\mathbb{R}^n$ which in addition
satisfy the first Bianchi identity. It is called the space of
algebraic curvature operators on $\mathbb{R}^n$. As a norm on this
space we use the classical Frobenius norm
$\|\R\|^2=\tr(\R^2)$. Similar constructions hold on the tangent bundle 
of a Riemannian manifold $(M,g)$ and give rise to the bundles
$S^2\Lambda^2TM$ and $S^2_B\Lambda^2TM$.

The curvature
tensor of a manifold $(M,g)$ will always be viewed as a section of the
bundle of curvture operators, $S^2_B\Lambda^2TM$. We follow the
convention of 
\cite{MR2415394} for the curvature operator, namely, the curvature
operator of a round sphere of radius 1 is the identity. 

We will use
$\R$, $\ric$ and $\scal$ to denote the curvature operator, Ricci
curvature and scalar curvature. When considering a Ricci flow
$(M,g(t))$, we will often not specify the dependence on $t$ of these
various curvature when no confusion is possible. We will write $\Id$
for the identity operator of $S^2_B\Lambda^2\mathbb{R}^n$ and $\id$
for the identity of $\mathbb{R}^n$.
 
Hamilton defined a bilinear map : 
\begin{align*}
  \#:S^2\Lambda^2\mathbb{R}^n\times S^2\Lambda^2\mathbb{R}^n &\to
  S^2\Lambda^2\mathbb{R}^n\\
  (\R,\Lop)&\mapsto \R\#\Lop
\end{align*}
whose expression can be found in \cite{MR862046} or \cite{MR2415394}.

If $g(t)$ is a
family of metric on $M$ evolving along the Ricci flow, Hamilton 
showed in \cite{MR862046} that in appropriate coordinates the
curvature operator $\R_{g(t)}$ satisfy the following evolution equation :
\[\ddt{\R}=\Delta \R + 2(\R^2+\R^\#),\]
where $\Delta$ is the connection laplacian and $\R^\#=\R\#\R$.

Removing the laplacian in this evolution equation, we obtain
Hamilton's ODE :
\[\dot{\R}=2(\R^2+\R^\#)=2Q(\R).\]

We will need the following algebraic fact about the $\#$ operator,
which was proved by Böhm and Wilking \cite[Lemma 2.1]{MR2415394} :
\begin{prop}\label{BWid}
  $\R+\R\#\Id=\ric\wedge\id$
\end{prop}
Here $\ric\wedge\id$ is the curvature operator defined by, for any $u$
and $v$ in $\mathbb{R}^n$ :
\[\ric\wedge\id(u\wedge v)=\frac{1}{2}(\ric(u)\wedge
v+u\wedge\ric(v)),\]
where $\ric$ is viewed as an operator on $\mathbb{R}^n$. In
particular, if $(\lambda_i)_{1\leq i\leq n}$ are the eigenvalues of
$\ric$ then the eigenvalues of $\ric\wedge\id$ are
$(\frac{\lambda_i+\lambda_j}{2})_{1\leq i<j\leq n}$.
\section{Proof of Theorem \ref{lboundthm}.}
\label{sec:proofthm}
\begin{proof}[Proof of Theorem \ref{lboundthm}]
  According to our hypothesis, if we define a new section of the
bundle of curvature operators $\Lop$ by :
\[\Lop=\R+\eps(\ph(t)+t\alpha\scal)\Id,\] 
it is enough to find a positive smooth function $\ph$, a constant $\alpha$ and a
time $T>0$, all depending only on $A$ and $B$ such that
$\Lop\in\mathcal{C}$ for $t\in [0,T)$. The fact that
$t\scal$ and $\ph$ are unformly bounded on $[0,T]$ will then give
the required bound. To ensure that $\Lop\in\mathcal{C}$ at time $0$, we
impose that $\ph(0)=1$. Since such lower bounds are likely to get
worse with time, we will assume that $\ph'\geq 0$.

To prove that $\Lop$ remains in $\mathcal{C}$, we will apply Hamilton's 
maximum principle for tensors \cite{MR862046}, or more precisely a
variant of it called maximum principle with avoidance set proved by
Chow and Lu in \cite[Theorem 4]{MR2042930}. This variant allows us to
use our a priori estimate on the scalar curvature (which is not
implied by the ODE) in the study of the ODE associated to the PDE satisfied
by $\Lop$. 

 We will impose conditions on $\ph$ and $\alpha$
during the proof and verify that these conditions can be fullfilled at
the end of the proof.

We first compute the evolution of $\Lop$ :
\begin{align*}
    \ddt{\Lop} = & \Delta \R+2Q(\R)
  +\eps(\varphi'+\alpha \scal+t\alpha(\Delta \scal+2|\ric|^2))\Id\\
  = & \Delta \Lop +2Q(\R)
  +\eps(\varphi' +\alpha \scal+2t\alpha |\ric|^2)\Id\\  
  = & \Delta \Lop + 2N(\Lop).
\end{align*}

We now have to show that $\mathcal{C}$ is preserved by the
differential equation $\dot{\Lop}=2N(\Lop)$. That is, given
$\Lop\in\partial\mathcal{C}$, we need to show that
$N(\Lop)\in\mathcal{C}$. Since $\mathcal{C}$ is preserved by
Hamilton's ODE, we know that $Q(\Lop)\in\mathcal{C}$ and we just need
to show (since $\mathcal{C}$ is convex) that
$D(\Lop)=N(\Lop)-Q(\Lop)\in\mathcal{C}$. This idea comes from the work
of Böhm and Wilking in \cite{MR2415394}.

We will in fact prove that
$D(\Lop)$ is a nonnegative curvature operator, which will be enough
since $\mathcal{C}$ contains the cone of nonnegative curvature operator.

Using Böhm and Wilking identity (proposition \ref{BWid}), we have :
\begin{align*}
  Q(\Lop)= & Q(\R)+2\eps(\ph+t\alpha\scal)(\R+\R\# \Id) +\eps^2(\ph+t\alpha\scal)^2 Q(\Id)\\
         = & Q(\R) +2\eps(\ph+t\alpha\scal)(\ric\wedge \id)+(n-1)
           \eps^2(\ph+t\alpha\scal)^2 \Id.
\end{align*}

We then compute $D(\Lop)$ :
\begin{align*}
  D(\Lop)= & N(\Lop)-Q(\Lop)\\
         = & \frac{\eps}{2}(\varphi' +\alpha \scal+2t\alpha
         |\ric|^2)\Id\\
           & -2\eps(\ph+t\alpha\scal)(\ric\wedge\id)-(n-1)
           \eps^2(\ph+t\alpha\scal)^2 \Id.
\end{align*}

In order to estimate the $2\ric\wedge \id$ term, we use that
$\Lop\in\mathcal{C}$ has nonnegative Ricci curvature, which gives that
$\ric\geq -(n-1)\eps(\ph+t\alpha\scal)\id$ as symmetric operators. Since $\tr(\ric)=\scal$, we
have :
\[-(n-1)\eps(\ph+t\alpha\scal)\id\leq\ric\leq (\scal+(n-1)^2\eps(\ph+t\alpha\scal))\id.\]
This implies :
\[2\ric\wedge \id\leq(\scal+(n-1)^2\eps(\ph+t\alpha\scal))\Id. \]
We now assume that :
\[\ph+t\alpha\scal\geq 0\qquad\text{condition (C1)}. \]

This allows us to estimate $D(\Lop)$ :
\begin{align*}
  D(\Lop)\geq &\ \frac{\eps}{2}(\varphi' +\alpha \scal)\Id\\
           & -\eps(\ph+t\alpha\scal)(\scal+(n-1)^2\eps(\ph+t\alpha\scal))\Id\\
           & -(n-1)\eps^2(\ph+t\alpha\scal)^2 \Id.
\end{align*}
We rearrange the
terms in the following way\footnote{We will drop the $\Id$'s in the next
  inequalities, here a real number $\alpha$ should be viewed as the
  operator $\alpha \Id$.} :
\begin{align*}
  D(\Lop)\geq&\ \frac{\eps}{2} \ph'\\
             & +\eps\scal((\frac{1}{2}-t\scal)\alpha-\ph)\\
             & -(2n-1)(n-1)\eps^2(\ph+t\alpha\scal)^2
\end{align*}
We now assume that :
\[0\leq (\frac{1}{2}-t\scal)\alpha-\ph \leq 1\qquad\text{condition (C2)},\]

Since $\scal\geq -\eps n(n-1)$ at $t=0$, it remains so as long as the
solution exists. Therefore we have :
\begin{align*}
  D(\Lop)\geq&\ \frac{\eps}{2}\ph'\\
             & -\eps^2 n(n-1)\\
             & -(2n-1)(n-1)\eps^2(\ph+t\alpha\scal)^2,
\end{align*}
and since $\eps\in [0,1]$ :
\begin{align*}
  \frac{1}{\eps^2}D(\Lop)\geq&\ \frac{\ph'}{2}\\
             & -n(n-1)\\
             & -(2n-1)(n-1)(\ph+t\alpha\scal)^2.
\end{align*}
We now use that $|t\scal|\leq A+Bt$ to get :
\begin{align*}
  \frac{1}{\eps^2}D(\Lop)\geq&\ \frac{\ph'}{2}\\
             & -n(n-1)\\
             & -(n-1)(2n-1)(\ph+\alpha(A+Bt))^2.
\end{align*}
To ensure that $D(\Lop)$ is a nonnegative operator, it is then enough to
show that :
\begin{equation*}
  \frac{\ph'}{2}
             -n(n-1)
             -(n-1)(2n-1)(\ph+\alpha(A+Bt))^2\geq
             0\qquad\text{condition (C3)}.
\end{equation*}

 We now have to find $\ph$, $\alpha$ and $T$
such that conditions (C1), (C2) and (C3) are satisfied on $[0,T]$.

Using again that $-n(n-1)t\leq t\scal \leq A+Bt$, we have that 
conditions (C1) and (C2) are implied by the following inequalities
which involves only $A$, $B$ and the dimension $n$ :
\begin{align*}
  \left .
  \begin{array}{l}
   (\frac{1}{2}-(A+Bt))\alpha-\ph \geq 0\\
   (\frac{1}{2}+tn(n-1))\alpha-\ph\leq 1\\ 
   \ph -n(n-1)t\alpha\geq 0
 \end{array} \right\}
   \qquad\text{condition (C4)}
\end{align*}
Looking at conditions (C4) at $t=0$, we see that it is fulfilled if
$\alpha$ belongs to $[\frac{2}{1-2A},4]$. We now impose that
$A<\frac{1}{4}$. Let $\alpha\in
(\frac{2}{1-2A},4)$, and $\ph(t)=1+\beta t$. Conditition (C4) is then
satisfied at time $0$ with strict inequalities.

We now choose $\beta$ big enough such that condition (C3) is fulfilled
with a strict inequality. By
continuity of $\varphi$, these conditions are still
fulfilled for $t$ in some small time interval $[0,T)$.

Our choices of $\ph$, $\alpha$ and $T$ depend only on $A$, $B$ and
$n$, the theorem is then proved.
\end{proof}
\section{First applications.}
\label{sec:applications}
\subsection{Gromov-Hausdorff converging sequences whose
  $\mathcal{C}$-curvature is bounded from below.}
\label{sec:grom-hausd-conv}
In this section, we prove Theorem \ref{GHrfthm}. We first state a
lemma which is of independent 
interest, the idea of using pseudolocality and convergence of the
isoperimetric profiles in the proof of the
following lemma was suggested to the author by
Gilles Carron :
\begin{lem}\label{GHvsRF}
 Let $(M_k,g_k)_{k\in\mathbb{N}}$ be a sequence of smooth compact
 $n$-dimensional Riemannian manifold which satisfies $\ric(g_k)\geq
 -(n-1)g_k$ and which GH-converges to a smooth compact $n$-dimensional
 Riemannian manifold $(M,g)$. 

 Then for every every $A>0$, there exist
 $k_0\in \mathbb{N}$, $B>0$ and $T>0$ such that, for any $k\geq k_0$ the Ricci flows
 $(M_k,g_k(t))$ whith initial condition $(M_k,g_k)$ exist at least on
 $[0,T)$ and satisfy :
 \begin{enumerate}
 \item $\|\R(g_k(t))\|\leq A/t+B$ for all $t\in(0,T)$,
 \item $\vol(B_{g_k(t)}(x,\sqrt{t}))\geq ct^{n/2}$ for all $t\in (0,T)$
   and $x\in M_i$.
 \end{enumerate}
 In particular, the Ricci flows $(M_k,g_k(t))_{t\in(0,T)}$ form a
 precompact sequence in the sense of Cheeger Gromov and Hamilton.
\end{lem}
\begin{proof}
  We want to apply Perelman's pseudolocality (\cite[Section
  10]{PerelEnt}, \cite[Theorem 30.1, Corollary 35.1]{KleiLott}) to get
  the two estimates of the lemma. The precompactness statement then
  follows from Hamilton's compactness theorem \cite{HamComp}.

  Let $A>0$ be fixed. We already know that for any $x\in M_k$,
  $\scal_{g_k}(x)\geq -n(n-1)$. Thus we just need to find some
  $r_0\in\left (0,(n(n-1))^{-1/2}\right ]$ such that any smooth domain
  $\Omega$ contained in a ball of radius $r_0$ in $M_k$ for $k$ large
  enough satisfies the almost Euclidean isoperimetric estimate :
  \begin{equation}
  |\partial\Omega|^{\frac{n}{n-1}}\geq (1-\delta)\gamma_n |\Omega|\label{eq:1}  
  \end{equation}
  where $\gamma_n$ is the euclidean isoperimetric constant and $\delta$
  is given by the pseudolocality theorem.

  To obtain this estimate, we will consider the isoperimetric profiles of
  the $(M_k,g_k)$'s, that will be denoted by $h_k(\beta)$. Since
  $(M,g)$ is smooth, by a result of Bérard and Meyer
  \cite[Appendice C]{MR690651},
  its isoperimetric profile $h(\beta)$ is equivalent to the
  euclidean one as $\beta$ goes to zero. Thus we can find, for any
  given $\varepsilon>0$, some $\rho>0$ such that :
  \[ \beta<\rho \Rightarrow\
  h(\beta)\geq (1-\varepsilon)
  \frac{\gamma_n}{\vol(M,g)^{\frac{1}{n}}} \beta^{\frac{n-1}{n}}\]

  We then use a result from Bayle thesis \cite{Bayle} : under non collapsing
  GH-convergence to a smooth manifold with Ricci curvature bounded
  from below, the ratio of the isoperimetric profiles $h_k/h$
  is going to $1$ uniformly on $(0,1)$. Then, for $i$ large enough :
  \[h_k\geq (1-\varepsilon)h.\]

  Let $\Omega\subset M_i$ be a smooth domain whose volume is less than
  $\rho \vol(M_k,g_k)$. We then have :
  \begin{align*}
  |\partial\Omega| & \geq 
  \vol(M,g)\times h_k\left(\frac{|\Omega|}{\vol(M,g)}\right)\\ 
  & \geq \vol(M,g)\times (1-\varepsilon)h
  \left(\frac{|\Omega|}{\vol(M,g)}\right)\\
  & \geq \vol(M,g)\times (1-\varepsilon)^2 
  \frac{\gamma_n}{\vol(M,g)^{\frac{1}{n}}}
  \left(\frac{|\Omega|}{\vol(M,g)}\right)^{\frac{n-1}{n}}\\
  & = (1-\varepsilon)^2 \gamma_n|\Omega|^{\frac{n-1}{n}}
  \end{align*}
  If we take $\eps$ small enough,  we get estimate \eqref{eq:1} for
  domains of volume less then $\rho \vol(M_k,g_k)$. 

  Now, using Colding's theorem on the continuity of
  volume \cite{ColVol}, for $k$ large enough, 
  $\vol(M_k,g_k)\geq V/2$  where $V$ is the volume of
  $(M,g)$. In particular, our almost Euclidean
  isoperimetric inequality is valid for domains of volume less than
  $\rho V/2$. Since the Ricci curvature is bounded from below, Bishop
  Gromov inequality gives us that :
  \[\vol(B_{g_k}(x,r))\leq V_{-1}(r)\]
  where $V_{-1}(r)$ is the volume a radius $r$ ball in the
  $n$-dimensional hyperbolic space. This shows that our isoperimetric
  inequality is valid for domains included in balls of radius less
  than $r_0$ where $r_0$ is such that $V_{-1}(r_0)=\rho V/2$.

  Finally, pseudolocality applies and we get the required
  bounds. 
\end{proof}
We now proove Theorem \ref{GHrfthm}.
\begin{proof}[Proof of Theorem \ref{GHrfthm}]
  We now consider a sequence $(M_k^n,g_k)$ of smooth compact manifolds
  whose $\mathcal{C}$-curvature is bounded from below by $-\Id$ and
  which in addition satisfy the assumptions of Lemma
  \ref{GHvsRF}. 

  Thanks to the previous lemma, we can find $i_0\in\mathbb{N}$, $T>0$ and
  a constant $B$ such that, for $k\geq k_0$, the Ricci flows
  $(M_k,g_k(t))$ satisfying $g_k(0)=g_k$ satisfy :
  \begin{equation}
  |\scal(g_k(t))|\leq\frac{1}{8t}+B\text{ for }t\in(0,T).\label{eq:3}
  \end{equation}

  We now use Theorem \ref{lboundthm} and the fact that $(M_k,g_k(0))$
  has $\mathcal{C}$-curvature bounded from below by $-\Id$ to find
  $T'>0$ and $K>0$ such that, for $t\in(0,T')$,
  \begin{equation}
  \R(g_k(t))\geq_\mathcal{C} -K\Id.\label{eq:2}  
  \end{equation}

  Since this implies that the
  Ricci curvature of $(M_k,g_k(t)))$ is bounded from below by
  $-(n-1)K$ on $[0,T')$, we can apply Lemma 6.1 in
  \cite{MS2009}. We get, for some constant $c>0$, that for $k\geq
  k_0$, $x,y\in M_k$ and $0<s\leq t <T'$ :
  \begin{equation}
  d_{g_k(s)}(x,y)-c(\sqrt{t}-\sqrt{s})\leq d_{g_k(t)}\leq
  e^{c(t-s)}d_{g_k(s)} \label{eq:4}
  \end{equation}
  where $d_{g_k(t)}$ is the distance function of $(M_k,g_k(t))$.

  Consider now a subsequence of the sequence
  $(M_k,g_k(t))_{t\in(0,T')}$ which converges in the sense of
  Cheeger-Gromov-Hamilton to a Ricci flow
  $(\tilde{M},\tilde{g}(t))_{t\in(0,T')}$. This flow also satisfies
  estimates \eqref{eq:3}, \eqref{eq:2} and \eqref{eq:4}.

  As in the proof of Theorem 9.2 in \cite{MS2009}, we can
  prove that the distances $d_{\tilde{g}(t)}$ uniformly converge as
  $t$ goes to zero to some distance $\tilde{d}$, which define the usual
  manifold topology on $\tilde{M}$, and that $(\tilde{M},\tilde{d})$
  is isometric to the GH-limit $(M,g)$ of the sequence
  $(M_k,g_k)$. In particular, $M$ and $\tilde{M}$ are homeomorphic. 
\end{proof}

\subsection{Manifolds with almost nonnegative $\mathcal{C}$-curvature.}
\label{sec:manif-with-almost}
We now proove Theorem \ref{almostNNCC}.
\begin{proof}
  By contradiction, take a sequence of counterexamples
  $(M_k,g_k)$ satisfying $\R\geq_\mathcal{C}-\eps_k\Id$, where
  $\eps_k$ goes to $0$,
  and the required bounds on the diameter and injectivity
  radius. We assume that none of the $M_k$ admits a metric with
  nonnegative $\mathcal{C}$-curvature. Without loss of generality, we
  assume that $\eps_k\leq 1$. 
  
  Since the injectivity radius and the Ricci curvature are bounded
  from below, we can use Anderson-Cheeger theorem
  \cite[Theorem 0.3]{AndersonCheeger2}. It gives us, for any $\eps>0$,
  some $r>0$ such that every ball $B$ of radius less than $r_0$ admit
  an harmonic coordinate chart $\ph_B:B\to\mathbb{R}^n$ with :
  \[\frac{1}{1+\eps}\ph_B^*\delta\leq g\leq (1+\eps)\ph_B^*\delta\]
  on $B$, where $\delta$ is the euclidean metric on $\mathbb{R}^n$.
  
  If we choose $\eps$ small enough, this control will give us an
  almost Euclidean isoperimetric estimate on balls of radius less than
  $r_0$. 

  Consider the sequence of Ricci flows $(M_k,g_k(t))$ such that
  $g_k(0)=g_k$. Pseudolocality gives : 
  \begin{itemize}
  \item each $(M_k,g_k(t))$ exists at least on $[0,T)$ where $T$ does
    not depend on $k$,
  \item for $t\in(0,T)$, $|\scal(g_k(t))|\leq \frac{1}{8t}+B$, where
    $B$ does not depend on $k$.
  \item the Ricci flows $(M_k,g_k(t))_{t\in(0,T)}$ form a precompact
    sequence in the sense of Cheeger-Gromov-Hamilton.
  \end{itemize}
  We can then apply Theorem \ref{lboundthm} to have that on some time
  interval $(0,T')\subset(0,T)$, $\R(g_k(t))\geq_\mathcal{C}
  -K\eps_k\Id$.

  Let $(M,g(t))_{t\in(0,T')}$ be the limit of a convergent subsequence
  of $(M_k,g_k(t))_{t\in(0,T')}$, it satisfies
  $\R(g(t))\geq_\mathcal{C} 0$ for $t\in(0,T')$. Now, since the Ricci
  curvature is bounded from below in time and along the sequence, we
  can find some constant $C$ such that :
  \[\diam(M_k,g_k(t))\leq e^{Ct}\diam(M_k,g_k(t))\leq e^{Ct}D\]
  for all $k\in\mathbb{N}$ and $t\in(0,T')$.
 
  This implies that $M$ is compact. Hence, we have a subsequence of
  $(M_k,g_k)$ all of whose elements are diffeomorphic to $M$, in
  particular, these elements admit a metric with non-negative
  $\mathcal{C}$-curvature. This is a contradiction.
\end{proof}

\section{Stronger results when operators in $\mathcal{C}$ have
  nonnegative sectional curvature}
\label{sec:manif-with-mathc}
{\bf We now assume that $\mathcal{C}$ contains the cone of curvature
  operators whose sectional curvature is nonnegative.}

As in the previous proofs, the crucial point is to
get an $A/t+B$ bound on the scalar curvature. We first state a lemma
which gives this bound when one has almost euclidean volume and
$\mathcal{C}$-curvature bounded from below at the initial time. This
lemma is a stronger version of Proposition 5.5 in \cite{2011arXiv1107.0606C}.
\begin{lem}\label{lemalmosteuc}
  For any dimension $n$, any $A\in(0,A_0(n))$, there exists $\kappa>0$,
  $\delta>0$, $\tilde{\kappa}>0$ and $T>0$ such that if $(M^n,g)$ is a compact Riemannian manifold such that :
  \begin{itemize}
  \item $\R\geq_\mathcal{C}-\kappa\Id$,
  \item $\forall x\in M^n,\quad \vol_g(B_g(x,1))\geq (1-\delta)\omega_n$,
  \end{itemize}
  where $\omega_n$ is the volume of the unit ball in $\mathbb{R}^n$.
  
  Then the Ricci flow $(M^n,g(t))$ with initial condition $(M^n,g)$
  exists at least on $[0,T)$ and satisfies :
  \begin{itemize}
  \item $\forall t\in(0,T)\quad \|\R(g(t))\|\leq \frac{A}{t}$,
  \item $\forall t\in(0,T)\quad \R(g(t))\geq-\tilde{\kappa}\Id$.
\end{itemize}

\end{lem}
\begin{proof}
  Since $|\scal|\leq 2\sqrt{n}\|\R\|$, we set
  $A_0(n)=\frac{1}{8\sqrt{n}}$. This ensures that if $\|\R\|\leq
  \frac{A}{t}$ we have the right estimate on the scalar curvature to
  apply Theorem \ref{lboundthm}. 

  The proof goes by contradiction. Fix $n\in\mathbb{N}$ and
  $A<A_0(n)$. Assume we can find a sequence of manifolds
  $(M_i,g_i)_{i\in\mathbb{N}}$ such that :
  \begin{itemize}
  \item $\R(g_i)\geq_\mathcal{C}-\delta_i\Id$,
  \item $\forall x\in M,\quad \vol_{g_i}(B_{g_i}(x,1))\geq
    (1-\delta_i)\omega_n$, 
  \end{itemize}
  for some sequence $(\delta_i)_{i\in\mathbb{N}}$ going to zero. And
  assume furthermore that the sequence $(t_i)_{i\in\mathbb{N}}$
  defined by $t_i=\sup\{\ t>0\ |\ \forall s\leq t,\ s\|R(g_i(s))\|\leq
  A\ \}$ goes to zero. Taking $i$ large enough, we can assume that
  $\delta_i\leq 1$ and $t_i$ is less than the time $T$ 
  given by Theorem \ref{lboundthm}, this ensures that for all
  $t\in[0,t_i]$ :
  \[\R(g_i(t))\geq_\mathcal{C} -K\delta_i \Id.\]
  
  With this lower bound, we can repeat word for word the proof of
  Proposition 5.5 in \cite{2011arXiv1107.0606C} and get the
  $\frac{A}{t}$ bound on the norm of the curvature operator on some
  time interval. The lower bound on $\mathcal{C}$-curvature is now
  given by Theorem \ref{lboundthm}.
\end{proof}
We now proove Theorem \ref{GHRFnnsec}. The Ricci flow of $(X,d)$ is
constructed as limit of the Ricci flows of the $(M_i,g_i)$, as in
\cite{MS2009}. 
\begin{proof}[Proof (of Theorem \ref{GHRFnnsec}).]
  Fix $\kappa$ and $\delta$ such that Lemma \ref{lemalmosteuc} apply
  with $A=\frac{1}{16\sqrt{n}}$. Consider a sequence $(M_i,g_i)$
  satifying the assumption of the 
  theorem. Using Lemma \ref{lemalmosteuc}, we have that the Ricci
  flows $(M_i,g_i(t))$ exist at least on $[0,T)$ and satisfy, for any
  $t$ in $[0,T)$ :
  \begin{itemize}
  \item $\|\R(g_i(t))\|\leq\frac{1}{16\sqrt{n}t}$,
  \item $\R(g_i(t))\geq_\mathcal{C}-K\Id$.
  \end{itemize}
  In addition, at time $t=0$, we have that any unit ball in any of the
  $M_i$'s has volume at least $(1-\delta)\omega_n$. This allow us to
  apply Lemma 6.1 and Corolarry 6.2 in \cite{MS2009} to get that, on
  some possibly smaller time interval $[0,T']$, we have the estimates,
  for some constant $C>0$ depending only on $\kappa$ and $\delta$ :
  \begin{itemize}
  \item $\forall x\in M_i,\
    \vol_{g_i}(B_{g_i}(x,1))\geq\frac{(1-\delta)\omega_n}{2}$,
  \item for $0<s\leq t\leq T'$, $d_{g_i(s)}-C(\sqrt{t}-\sqrt{s})\leq
    d_{g_i(t)}\leq e^{C(t-s)}d_{g_i(s)}$,
  \end{itemize}
  where $d_{g_i(t)}$ is the distance on $M_i$ induced by the metric $g_i(t)$.

  We then argue as in the proof of Theorem 9.2 in \cite{MS2009} and
  get that the sequence of Ricci flows $(M_i,g_i(t))_{t\in(0,T')}$ has
  a convergent subsequence whose limit $(M,g(t))_{t\in(0,T')}$ is a
  Ricci flow of the Gromov-Hausdorff limit $(X,d)$ of the sequence
  $(M_i,g_i)$ in the sense that it satisfies the conclusions of
  Theorem \ref{GHRFnnsec}. 
\end{proof}

We now go on with Theorem \ref{almostnnccnnsc}.
\begin{proof}[Proof of Theorem \ref{almostnnccnnsc}.]
  Let $\delta>0$ and $\kappa>0$ be the constants given by Lemma
  \ref{lemalmosteuc} with 
  $A=\frac{1}{16\sqrt{n}}$. Fix $D>0$.

  As in the proof of Theorem \ref{almostNNCC}, consider a sequence of
  manifolds $(M_i,g_i)$ with :
  \begin{itemize}
  \item $\forall x\in M_i,\
    \vol_{g_i}(B_{g_i}(x,1))\geq(1-\delta)\omega_n$,
  \item $\R\geq_\mathcal{C}\ -\eps_i\Id$,
  \item $\diam(M_i,g_i)\leq D$,
 \end{itemize}
  where $\eps_i$ goes to $0$ as $i$ goes to infinity. Assume
  furthermore that none of the $M_i$ admits a metric with nonnegative
  $\mathcal{C}$-curvature. Without loss of
  generality, we can assume that $\eps_i\leq\min(\kappa,1)$.

  Arguing as in the proofs of Theorem \ref{GHRFnnsec} and Theorem
  \ref{almostNNCC}, we get that the 
  Ricci flows $(M_i,g_i(t))$ starting at $(M_i,g_i)$ exist at least on
  $[0,T)$, form a precompact family in the sense of
  Cheeger-Gromov-Hamilton, and satisfy :
  \begin{itemize}
  \item $\R\geq_\mathcal{C} -K\eps_i\Id$,
  \item $\diam(M_i,g_i(t))\leq e^{Ct}D$.
  \end{itemize}
  In particular any limit of a convergent subsequence will be compact
  and have nonnegative $\mathcal{C}$-curvature. Thus the sequence
  contains manifolds which admit metrics with
  $\mathcal{C}$-nonnegative curvature. This is a contradiction.  
\end{proof}

{\sc Institut Fourier, 100 rue des maths, 38402 St Martin d'Hères.}

{\it Email adress:} {\tt thomas.richard@ujf-grenoble.fr}
\end{document}